\documentclass[10pt,a4paper,twoside,leqno]{article}
\usepackage{amsfonts,amssymb}
\usepackage{eufrak}
\usepackage{amscd}

\newtheorem{defn}{Definition}

\vspace{0.5cm}

 4

\font\erm=cmr8

\usepackage{graphicx}

\begin{document}

\thispagestyle{empty}

\noindent {\bf New type Stirling like numbers - an email  style
letter }

\vspace{0.3cm} {\it A. Krzysztof Kwa\'sniewski}

\vspace{0.1cm}

{\erm High School of Mathematics and Applied Informatics}

{\erm  Kamienna 17, PL-15-021 Bia\l ystok, Poland}

{\erm e-mail: kwandr@wp.pl

\vspace{0.1cm}

AMS Classification Numbers: 05C20, 11C08, 17B56 .

\vspace{0.1cm}

\noindent Key Words: prefab, cobweb poset, Whitney numbers,
Fibonacci like sequences

\vspace{0.1cm}

\noindent presentation  (November $2005$) at the Gian-Carlo Rota
Polish Seminar\\
\noindent \emph{http://ii.uwb.edu.pl/akk/sem/sem\_rota.htm}

\vspace{0.2cm}

\noindent The notion of the Fibonacci cobweb poset from [1] has
been naturally extended to any admissible sequence $F$ in [2]
where it was also recognized that the celebrated  prefab notion of
Bender and Goldman [3] - (see also [4,5]) - admits such an
extension so as to encompass the new type combinatorial objects
from [2] as leading examples. Recently the present author  had
introduced also [6] two natural partial orders in there: one
$\leq$ in grading-natural subsets of cobweb`s prefabs sets [2] and
in the second proposal one endows the set sums of the so called
"prefabiants" with such another partial order that one arrives at
Bell-like numbers including Fibonacci triad sequences introduced
by the present author in [7]. Here we quote the basic observations
concerning the new type Stirling like numbers as they appear in
[6]. For more on notation, Stirling like numbers of the first kind
and for proofs - see [6].

\vspace{1mm} \noindent\textbf{The overall $F$-independent class of
p.o. set structure.} Let the family $S$ of combinatorial objects
($prefabiants$) consists of all layers $\langle\Phi_k \rightarrow
\Phi_n \rangle,\quad k<n,\quad k,n \in N\cup\{0\}\equiv Z_\geq$
and an empty prefabiant $i$. The set $\wp$ of prime objects
consists of all sub-posets $\langle\Phi_0 \rightarrow \Phi_m
\rangle$  i.e. all $P_m$`s $m \in N\cup\{0\}\equiv Z_\geq$
constitute from now on a family of prime $prefabiants$ [2]. Layer
is considered here to be the set of  all max-disjoint isomorphic
copies (iso-copies) of $P_m = P_{n-k}$ [2]. Consider then now  the
partially ordered family $S$ of these layers considered to be sets
of  all max-disjoint isomorphic copies (iso-copies) of  prime
prefabiants $P_{m}= P_{n-k}$. For any $F$-sequence determining
cobweb poset [2] let us define in $S$ \textbf{the same} partial
order relation as follows.
\begin{defn}
$$
\langle \Phi_k \rightarrow \Phi_n \rangle \leq
\langle\Phi_{k^\clubsuit} \rightarrow \Phi_{n^\clubsuit}\rangle
\quad \equiv\quad k \leq k{^\clubsuit} \quad\wedge\quad n\leq
n^{\clubsuit}.
$$
\end{defn}
For convenience reasons we shall also adopt and use the following
notation: $$ \langle\Phi_k \rightarrow \Phi_n \rangle = p_{k,n}.
$$
In what follows we shall consider  the subposet  $\langle\ P_{k,n}
, \leq \rangle$ where $ P_{k,n}= [p_{o,o} ,p_{k,n}].$ Then
according to [6] we observe  the following.

\noindent {\bf Observation 1.} The size $|P_{k,n}|$ of $P_{k,n}$ $
=|\{\langle l,m \rangle ,\quad 0 \leq l \leq k \quad\wedge \quad 0
\leq m \leq n \quad\wedge \quad k\leq n \}| = (n-k)(k + 1) +
\frac{k(k+1)}{2}.$

\noindent {\bf Observation 2.} The number of maximal chains in
$\langle\ P_{k,n} , \leq \rangle$ is equal to the number $d(k,n)$
of $0$ - dominated strings of binary sequences
$$d(k,n) = \frac{n+1-k}{n} \left( \begin{array}{c}
{k+n}\\n\end{array}\right). $$ Recall that $( d(k,n) )$ infinite
matrix`s  diagonal elements are equal to the \textbf{Catalan}
numbers $C(n)$. The poset $\langle\ P_{k,n} , \leq \rangle$ is
naturally graded. $\langle\ P_{k,n} , \leq \rangle$ poset`s
maximal chains are all of equal size (Dedekind -Jordan property)
therefore the rang function is defined.

\vspace{1mm}

\noindent {\bf Observation 3.} The rang $r(P_{k,n})$  of $P_{k,n}$
= number of elements in maximal chains $P_{k,n}$ minus one
$=k+n-1$. \noindent The rang $r(p_{l,m})$ of $\pi = p_{l,m} \in
P_{k,n}$ is defined accordingly: $r(p_{l,m})= l+m-1.$

\vspace{1mm}

\noindent Accordingly Whitney numbers $W_k(P_{l,m})$  of the
second kind are defined as follows (association: $ n
\leftrightarrow  \langle l,m \rangle $)

\begin{defn}
$$
W_k(P_{l,m})= \sum_{\pi\in P_{l,m}, r(\pi)=k} 1 \quad\equiv\quad
S(k,\langle l,m \rangle).
$$
\end{defn}
We shall identify  $W_k(P_{l,m})$ with $S(k,\langle l,m \rangle)$
called and  viewed at  as\textbf{ Stirling - like numbers} of the
second kind  of the naturally graded poset $\langle\ P_{k,n} ,
\leq \rangle$. Let us define also the corresponding Bell-like
numbers $ B(\langle l,m \rangle)$ of the naturally graded poset
$\langle\ P_{k,n} , \leq \rangle$.

\begin{defn}
$$
B(\langle l,m \rangle)= \sum_{k=0}^{l+m}S(k,\langle l,m \rangle).
$$
\end{defn}
{\bf Observation 4.}
$$ B(\langle l,m \rangle)=
|P_{l,m}| = \frac{k(k+1)}{2} + (n-k)(k+1).
$$
\textbf{The $F$-dependent, $F$-labelled  class of p.o. set
structures.} Let us consider now  prefabiants` set sums with an
appropriate another partial order so as to arrive at Bell-like
numbers including Fibonacci triad sequences introduced recently by
the present author in [7]. \noindent Let $F$ be any \textit{
"GCD-morphic"} sequence [2]. This means that $GCD[F_n,F_m] =
F_{GCD[n,m]}$ where $GCD$ stays for Greatest Common Divisor
mapping. We define the $F$-\textit{dependent} finite partial
ordered set $P(n,F)$ as the set of \textbf{prime} prefabiants
$P_l$ given by the sum below.

\begin{defn}
$$ P(n,F) =\bigcup_{0\leq p}\langle \Phi_p \rightarrow \Phi_{n-p} \rangle = \bigcup_{0\leq l}P_{n-l}$$
\end{defn}
with the partial order relation defined  for $n-2l\leq 0$
according to
\begin{defn}
$$
P_l \leq P_{\hat l} \quad \equiv \quad l\leq \hat l,\quad P_{\hat
l}, P_l\in \langle \Phi_l \rightarrow \Phi_{n-l} \rangle.
$$
\end{defn}
Recall that \textbf{rang of} $P_l$ \textbf{is} $l$. Note that
$\langle \Phi_l \rightarrow \Phi_{n-l} \rangle = \emptyset$ for
$n-2l \leq 0$. The Whitney numbers of the second kind are
introduce accordingly.

\vspace{2mm}

\begin{defn}
$$
W_k(P_{n,F})= \sum_{\pi\in P(n,F), r(\pi)=k} 1 \equiv S(n,n-k,F).
$$
\end{defn}
Right from the definitions above we infer that:

\vspace{1mm}

\noindent \textbf{Observation 5.}
$$
W_k(P_{n,F})= \sum_{\pi\in P(n,F), r(\pi)=k} 1 \equiv S(n,n-k,F)=
\left( \begin{array}{c} {n-k}\\k\end{array}\right)_{F}.
$$
Referring to the classical examples from [8] we identify
$W_k(P(n,F))$  with $S(n,n-k,F)$ called the Stirling - like
numbers of the second kind of the $P$. $P$ by construction
displays self-similarity property with respect to its prime
prefabiants sub - posets $P_n = P(n,F)$. Consequently for any
$GCD$-morphic sequence $F$ (see: [2])  we define the corresponding
Bell-like numbers $ B_n(F)$ of the poset $P(n,F)$ as follows.
\begin{defn}
$$B_n(F)= \sum_{k\geq0}S(n,k,F).$$
\end{defn}
Due to the investigation in [7] we have right now at our disposal
all corresponding results of [7] as the following identification
with  special case of $\langle \alpha,\beta,\gamma\rangle$ -
Fibonacci sequence $\langle
F_n^{[\alpha,\beta,\gamma]}\rangle_{n\geq 0}$ defined in [7]
holds.

\vspace{1mm}

\noindent {\bf Observation 6.}
$$B_n(F)\equiv F_{n+1}^{[\alpha=0,\beta=0,\gamma=0].}$$
Proof: See the Definition 2.2. from [7].

\noindent\textbf{ Recurrence relations.} Recurrence relations for
$\langle \alpha,\beta,\gamma\rangle$ - Fibonacci sequences
$F_n^{[\alpha,\beta=,\gamma]}$ are to be found in [7] - formula
(9). Compare also with the special case formula (7) in [9].

\vspace{1mm}

\noindent\ \textbf{Remark.} As seen from the identification
Observation 6. the special cases of $\langle
\alpha,\beta,\gamma\rangle$ - Fibonacci sequences
$F_n^{[\alpha,\beta,\gamma]}$  gain \textbf{additional} with
respect to [7] combinatorial interpretation in terms Bell-like
numbers as just sums of Whitney numbers of the poset $P(n,F)$.
This adjective \textit{"additional"} applies spectacularly to
Newton binomial connection constants between bases
$\langle(x-1)^k\rangle_{k\geq0}$ and $\langle x^n\rangle_{n\geq0}$
as these are Whitney numbers of the numbers from $[n]$ chain i.e.
Whitney numbers of the poset $\langle [n], \leq \rangle.$ For
other elementary "shining brightly"  examples see Joni , Rota and
Sagan excellent presentation in [8].

\begin
{thebibliography}{99}
\parskip 0pt

\bibitem{1}
A. K.  Kwa\'sniewski, {\it  Comments on  combinatorial
interpretation of fibonomial coefficients - an email  style
letter} Bulletin of the Institute of Combinatorics and its
Applications {\bf 42}  (2004), 10-11.

\bibitem{2}
A. K.  Kwa\'sniewski, {\it Cobweb posets as noncommutative
prefabs} submitted for publication ArXiv : math.CO/0503286 (2005)

\bibitem{3}
E. Bender, J. Goldman   {\it Enumerative uses of generating
functions} , Indiana Univ. Math.J. {\bf 20} 1971), 753-765.

\bibitem{4}
D. Foata and M. Sch"utzenberger, Th'eorie g'eometrique des
polynomes euleriens, (Lecture Notes in Math., No. 138).
Springer-Verlag, Berlin and New York, 1970.

\bibitem{5}
A. Nijenhuis and H. S. Wilf, Combinatorial Algorithms, 2nd ed.,
Academic Press, New York, 1978.

\bibitem{6}
A.K. Kwa\'sniewski, {\it Prefab posets` Whitney numbers} ArXiv:
math.CO/0510027 , 3 Oct 2005

\bibitem{7}
A. K.  Kwa\'sniewski, {\it Fibonacci-triad sequences} Advan. Stud.
Contemp. Math. {\bf 9} (2) (2004),109-118.

\bibitem{8}
S.A. Joni ,G. C. Rota, B. Sagan {\it From sets to functions: three
elementary examples} Discrete Mathematics {\bf 37} (1981),
193-2002.

\end{thebibliography}



\end{document}